\documentclass[11pt]{article}

\usepackage{latexsym}
\usepackage{amssymb}

\newtheorem{thm}{Theorem}[section]

\newtheorem{prop}{Proposition}[section]

\newtheorem{dfn}{Definition}[section]

\begin{document}
{

\begin{center}
\Large\bf
Devinatz's moment problem: a description of all solutions.
\end{center}
\begin{center}
\bf S.M. Zagorodnyuk
\end{center}

\section{Introduction.}
We shall study the following problem: to find a non-negative Borel measure $\mu$ in a strip
$$ \Pi = \{ (x,\varphi):\ x\in \mathbb{R},\ -\pi\leq \varphi < \pi \}, $$
such that
\begin{equation}
\label{f1_1}
\int_\Pi x^m e^{in\varphi} d\mu = s_{m,n},\qquad m\in \mathbb{Z}_+, n\in \mathbb{Z},
\end{equation}
where $\{ s_{m,n} \}_{m\in \mathbb{Z}_+, n\in \mathbb{Z}}$ is a given sequence of complex numbers.
We shall refer to this problem as to {\bf the Devinatz moment problem}.

\noindent
A.~Devinatz was the first who introduced and studied this moment problem~\cite{cit_1000_D}.
He obtained the necessary and sufficient conditions of solvability for the moment problem~(\ref{f1_1})
and gave a sufficient condition for the moment problem to be determinate~\cite[Theorem 4]{cit_1000_D}.

\noindent
Our aim here is threefold. Firstly, we present a  new proof of the Devinatz solvability criterion.
Secondly, we describe canonical solutions of the Devinatz moment problem (see the definition
below).
Finally, we describe all solutions of the Devinatz moment problem.
We shall use an abstract operator approach~\cite{cit_1500_Z} and results of Godi\v{c}, Lucenko and
Shtraus~\cite{cit_2000_GL},\cite[Theorem 1]{cit_3000_GP},\cite{cit_4000_S}.

{\bf Notations. } As usual, we denote by $\mathbb{R},\mathbb{C},\mathbb{N},\mathbb{Z},\mathbb{Z}_+$
the sets of real numbers, complex numbers, positive integers, integers and non-negative integers,
respectively.
For a subset $S$ of the complex plane we denote by $\mathfrak{B}(S)$ the set of all Borel subsets of $S$.
Everywhere in this paper, all Hilbert spaces are assumed to be separable. By
$(\cdot,\cdot)_H$ and $\| \cdot \|_H$ we denote the scalar product and the norm in a Hilbert space $H$,
respectively. The indices may be ommited in obvious cases.
For a set $M$ in $H$, by $\overline{M}$ we mean the closure of $M$ in the norm $\| \cdot \|_H$. For
$\{ x_k \}_{k\in T}$, $x_k\in H$, we write
$\mathop{\rm Lin}\nolimits \{ x_k \}_{k\in T}$ for the set of linear combinations of vectors $\{ x_k \}_{k\in T}$
and $\mathop{\rm span}\nolimits \{ x_k \}_{k\in T} =
\overline{ \mathop{\rm Lin}\nolimits \{ x_k \}_{k\in T} }$.
Here $T := \mathbb{Z}_+ \times \mathbb{Z}$, i.e. $T$ consists of pairs $(m,n)$,
$m\in \mathbb{Z}_+$, $n\in\mathbb{Z}$.
The identity operator in $H$ is denoted by $E$. For an arbitrary linear operator $A$ in $H$,
the operators $A^*$,$\overline{A}$,$A^{-1}$ mean its adjoint operator, its closure and its inverse
(if they exist). By $D(A)$ and $R(A)$ we mean the domain and the range of the operator $A$.
By $\sigma(A)$, $\rho(A)$ we denote the spectrum of $A$ and the resolvent set of $A$, respectively.
We denote by $R_z (A)$ the resolvent function of $A$, $z\in \rho(A)$.
The norm of a bounded operator $A$ is denoted by $\| A \|$.
By $P^H_{H_1} = P_{H_1}$ we mean the operator of orthogonal projection in $H$ on a subspace
$H_1$ in $H$. By $\mathbf{B}(H)$ we denote the set of all bounded operators in $H$.

\section{Solvability.}
Let a moment problem~(\ref{f1_1}) be given.
Suppose that the moment problem has a solution $\mu$. Choose an arbitrary power-trigonometric
polynomial $p(x,\varphi)$ of the following form:
\begin{equation}
\label{f1_2}
\sum_{m=0}^\infty \sum_{n=-\infty}^\infty \alpha_{m,n} x^m e^{in\varphi},\qquad \alpha_{m,n}\in \mathbb{C},
\end{equation}
where all but finite number of coefficients $\alpha_{m,n}$ are zeros.
We can write
$$ 0 \leq \int_\Pi |p(x,\varphi)|^2 d\mu =
\int_\Pi  \sum_{m=0}^\infty \sum_{n=-\infty}^\infty \alpha_{m,n} x^m e^{in\varphi}
\overline{
\sum_{k=0}^\infty \sum_{l=-\infty}^\infty \alpha_{k,l} x^k e^{il\varphi}
} d\mu $$
$$ = \sum_{m,n,k,l} \alpha_{m,n}\overline{\alpha_{k,l}} \int_\Pi x^{m+k} e^{i(n-l)\varphi} d\mu =
\sum_{m,n,k,l} \alpha_{m,n}\overline{\alpha_{k,l}} s_{m+k,n-l}. $$
Thus, for arbitrary complex numbers $\alpha_{m,n}$ (where all but finite numbers are zeros) we have
\begin{equation}
\label{f2_1}
\sum_{m,k=0}^\infty \sum_{n,l=-\infty}^\infty \alpha_{m,n}\overline{\alpha_{k,l}} s_{m+k,n-l} \geq 0.
\end{equation}
Let $T = \mathbb{Z}\times \mathbb{Z}_+$ and for $t,r\in T$, $t=(m,n)$, $r=(k,l)$, we set
\begin{equation}
\label{f2_2}
K(t,r) = K((m,n),(k,l)) = s_{m+k,n-l}.
\end{equation}
Thus, for arbitrary elements $t_1,t_2,...,t_n$ of $T$ and
arbitrary complex numbers $\alpha_1,\alpha_2,...,\alpha_n$, with $n\in \mathbb{N}$, the following inequality holds:
\begin{equation}
\label{f2_3}
\sum_{i,j=1}^n K(t_i,t_j) \alpha_{i} \overline{\alpha_j} \geq 0.
\end{equation}
The latter means that $K(t,r)$ is a positive matrix in the sense of E.H.~Moore \cite[p.344]{cit_5000_A}.

Suppose now that a Devinatz moment problem is given and conditions~(\ref{f2_1}) (or what is the same
conditions~(\ref{f2_3})) hold. Let us show that the moment problem has a solution.
We shall use the following important fact (e.g.~\cite[pp.361-363]{cit_6000_AG}).
\begin{thm}
\label{t2_1}
Let $K = K(t,r)$ be a positive matrix on $T=\mathbb{Z}\times \mathbb{Z}_+$.
Then there exist a separable Hilbert space $H$ with a scalar product $(\cdot,\cdot)$ and
a sequence $\{ x_t \}_{t\in T}$ in $H$, such that
\begin{equation}
\label{f2_4}
K(t,r) = (x_t,x_r),\qquad t,r\in T,
\end{equation}
and $\mathop{\rm span}\nolimits\{ x_t \}_{t\in T} = H$.
\end{thm}
{\bf Proof. }
Consider an arbitrary infinite-dimensional linear vector space $V$ (for example, we can choose a space of complex
sequences $(u_n)_{n\in \mathbb{N}}$, $u_n\in \mathbb{C}$).
Let $X = \{ x_t \}_{t\in T}$ be an arbitrary infinite sequence of linear independent elements
in $V$ which is indexed by elements  of $T$.
Set $L_X = \mathop{\rm Lin}\nolimits\{ x_t \}_{t\in T}$. Introduce the following functional:
\begin{equation}
\label{f2_5}
[x,y] = \sum_{t,r\in T} K(t,r) a_t\overline{b_r},
\end{equation}
for $x,y\in L_X$,
$$ x=\sum_{t\in T} a_t x_t,\quad y=\sum_{r\in T} b_r x_r,\quad a_t,b_r\in \mathbb{C}. $$
Here all but finite number of indices $a_t,b_r$ are zeros.

\noindent
The set $L_X$ with $[\cdot,\cdot]$ will be a pre-Hilbert space. Factorizing and making the completion
we obtain the  required space $H$ (\cite[p. 10-11]{cit_7000_B}).
$\Box$

By applying this theorem we get that there exist a Hilbert space $H$ and a sequence
$\{ x_{m,n} \}_{m\in \mathbb{Z}_+, n\in \mathbb{Z}}$, $x_{m,n}\in H$, such that
\begin{equation}
\label{f2_6}
(x_{m,n}, x_{k,l})_H = K((m,n),(k,l)),\qquad m,k\in \mathbb{Z}_+,\ n,l\in \mathbb{Z}.
\end{equation}
Set $L = \mathop{\rm Lin}\nolimits\{ x_{m,n} \}_{(m,n)\in T}$.
We introduce the following operators
\begin{equation}
\label{f2_7}
A_0 x = \sum_{(m,n)\in T} \alpha_{m,n} x_{m+1,n},
\end{equation}
\begin{equation}
\label{f2_8}
B_0 x = \sum_{(m,n)\in T} \alpha_{m,n} x_{m,n+1},
\end{equation}
where
\begin{equation}
\label{f2_9}
x = \sum_{(m,n)\in T} \alpha_{m,n} x_{m,n} \in L.
\end{equation}
We should show that these definitions are correct.
Indeed, suppose that the element $x$ in~(\ref{f2_9}) has another representation:
\begin{equation}
\label{f2_10}
x = \sum_{(k,l)\in T} \beta_{k,l} x_{k,l}.
\end{equation}
We can write
$$ \left( \sum_{(m,n)\in T} \alpha_{m,n} x_{m+1,n}, x_{a,b} \right) =
\sum_{(m,n)\in T} \alpha_{m,n} K((m+1,n),(a,b)) $$
$$= \sum_{(m,n)\in T} \alpha_{m,n} s_{m+1+a,n-b} = \sum_{(m,n)\in T} \alpha_{m,n} K((m,n),(a+1,b)) $$
$$ = \left(\sum_{(m,n)\in T} \alpha_{m,n} x_{m,n}, x_{a+1,b} \right) = (x,x_{a+1,b}), $$
for arbitrary $(a,b)\in T$.
In the same manner we get
$$ \left(\sum_{(k,l)\in T} \beta_{k,l} x_{k+1,l}, x_{a,b} \right) = (x,x_{a+1,b}). $$
Since $\mathop{\rm span}\nolimits\{ x_{a,b} \}_{(a,b)\in T} = H$, we get
$$ \sum_{(m,n)\in T} \alpha_{m,n} x_{m+1,n} =
\sum_{(k,l)\in T} \beta_{k,l} x_{k+1,l}. $$
Thus, the operator $A_0$ is defined correctly.

\noindent
We can write
$$ \left\| \sum_{(m,n)\in T} (\alpha_{m,n}-\beta_{m,n}) x_{m,n+1} \right\|^2 $$
$$= \left( \sum_{(m,n)\in T} (\alpha_{m,n}-\beta_{m,n}) x_{m,n+1},
\sum_{(k,l)\in T} (\alpha_{k,l}-\beta_{k,l}) x_{k,l+1} \right) $$
$$ = \sum_{(m,n),(k,l)\in T} (\alpha_{m,n}-\beta_{m,n}) \overline{(\alpha_{k,l}-\beta_{k,l})}
K((m,n+1),(k,l+1)) $$
$$= \sum_{(m,n),(k,l)\in T} (\alpha_{m,n}-\beta_{m,n}) \overline{(\alpha_{k,l}-\beta_{k,l})}
K((m,n),(k,l))  $$
$$= \left( \sum_{(m,n)\in T} (\alpha_{m,n}-\beta_{m,n}) x_{m,n},
\sum_{(k,l)\in T} (\alpha_{k,l}-\beta_{k,l}) x_{k,l} \right) = 0. $$
Consequently, the operator $B_0$ is defined correctly, as well.

Choose an arbitrary $y = \sum_{(a,b)\in T} \gamma_{a,b} x_{a,b} \in L$. We have
$$ (A_0 x,y) = \sum_{m,n,a,b} \alpha_{m,n}\gamma_{a,b} (x_{m+1,n},x_{a,b}) =
\sum_{m,n,a,b} \alpha_{m,n}\gamma_{a,b} K((m+1,n),(a,b)) $$
$$ = \sum_{m,n,a,b} \alpha_{m,n}\gamma_{a,b} K((m,n),(a+1,b)) =
 \sum_{m,n,a,b} \alpha_{m,n}\gamma_{a,b} (x_{m,n},x_{a+1,b}) =  (x,A_0 y). $$
Thus, $A_0$ is a symmetric operator. Its closure we denote by $A$.
On the other hand, we have
$$ (B_0 x,B_0 y) = \sum_{m,n,a,b} \alpha_{m,n}\overline{\gamma_{a,b}} (x_{m,n+1},x_{a,b+1}) =
\sum_{m,n,a,b} \alpha_{m,n}\overline{\gamma_{a,b}} K((m,n+1),(a,b+1)) $$
$$ = \sum_{m,n,a,b} \alpha_{m,n}\overline{\gamma_{a,b}} K((m,n),(a,b)) =
 \sum_{m,n,a,b} \alpha_{m,n}\overline{\gamma_{a,b}} (x_{m,n},x_{a,b}) =  (x,y). $$
In particular, this means that $B_0$ is bounded. By continuity we extend $B_0$ to a bounded
operator $B$ such that
$$ (Bx,By) = (x,y),\qquad x,y\in H. $$
Since $R(B_0)=L$ and $B_0$ has a bounded inverse, we have $R(B)=H$. Thus, $B$ is a unitary operator in $H$.

Notice that operators $A_0$ and $B_0$ commute. It is straightforward to check that $A$ and $B$ commute:
\begin{equation}
\label{f2_11}
AB x = BA x,\qquad x\in D(A).
\end{equation}
Consider the following operator:
\begin{equation}
\label{f2_12}
J_0 x = \sum_{(m,n)\in T} \overline{\alpha_{m,n}} x_{m,-n},
\end{equation}
where
\begin{equation}
\label{f2_13}
x = \sum_{(m,n)\in T} \alpha_{m,n} x_{m,n} \in L.
\end{equation}
Let us check that this definition is correct. Consider another representation for $x$ as in~(\ref{f2_10}).
Then
$$ \left\| \sum_{(m,n)\in T} (\overline{\alpha_{m,n}} - \overline{\beta_{m,n}}) x_{m,-n} \right\|^2 $$
$$= \left( \sum_{(m,n)\in T} \overline{ (\alpha_{m,n}-\beta_{m,n}) } x_{m,-n},
\sum_{(k,l)\in T} \overline{ (\alpha_{k,l}-\beta_{k,l}) } x_{k,-l} \right) $$
$$ = \sum_{(m,n),(k,l)\in T} \overline{(\alpha_{m,n}-\beta_{m,n})} (\alpha_{k,l}-\beta_{k,l})
K((m,-n),(k,-l)) $$
$$= \overline{ \sum_{(m,n),(k,l)\in T} (\alpha_{m,n}-\beta_{m,n}) \overline{(\alpha_{k,l}-\beta_{k,l})}
K((m,n),(k,l)) }  $$
$$= \overline{ \left( \sum_{(m,n)\in T} (\alpha_{m,n}-\beta_{m,n}) x_{m,n},
\sum_{(k,l)\in T} (\alpha_{k,l}-\beta_{k,l}) x_{k,l} \right) } = 0. $$
Thus, the definition of $J_0$ is correct.
For an arbitrary $y = \sum_{(a,b)\in T} \gamma_{a,b} x_{a,b} \in L$ we can write
$$ (J_0 x,J_0 y) = \sum_{m,n,a,b} \overline{\alpha_{m,n}}\gamma_{a,b} (x_{m,-n},x_{a,-b}) =
\sum_{m,n,a,b} \overline{\alpha_{m,n}}\gamma_{a,b} K((m,-n),(a,-b)) $$
$$ = \sum_{m,n,a,b} \overline{\alpha_{m,n}} \gamma_{a,b} K((a,b),(m,n)) =
 \sum_{m,n,a,b} \overline{\alpha_{m,n}}\gamma_{a,b} (x_{a,b},x_{m,n}) =  (y,x). $$
In particular, this implies that $J_0$ is bounded. By continuity we extend $J_0$ to a bounded antilinear
operator $J$ such that
$$ (Jx,Jy) = (y,x),\qquad x,y\in H. $$
Moreover, we get $J^2 = E_H$. Consequently, $J$ is a conjugation in $H$ (\cite{cit_8000_S}).

\noindent
Notice that $J_0$ commutes with $A_0$. It is easy to check that
\begin{equation}
\label{f2_14}
AJ x = JA x,\qquad x\in D(A).
\end{equation}
On the other hand, we have $J_0 B_0 = B_0^{-1} J_0$. By continuity we get
\begin{equation}
\label{f2_15}
JB = B^{-1}J.
\end{equation}
Consider the Cayley transformation of the operator A:
\begin{equation}
\label{f2_16}
V_A := (A+iE_H)(A-iE_H)^{-1},
\end{equation}
and set
\begin{equation}
\label{f2_17}
H_1 := \Delta_A(i),\ H_2 := H\ominus H_1,\ H_3:= \Delta_A(-i),\ H_4 := H\ominus H_3.
\end{equation}
\begin{prop}
\label{p2_1}
The operator $B$ reduces subspaces $H_i$, $1\leq i\leq 4$:
\begin{equation}
\label{f2_18}
BH_i = H_i,\qquad 1\leq i\leq 4.
\end{equation}
Moreover, the following equality holds:
\begin{equation}
\label{f2_19}
BV_Ax = V_ABx,\qquad x\in H_1.
\end{equation}
\end{prop}
{\bf Proof. } Choose an arbitrary $x\in \Delta_A(z)$, $x=(A-zE_H)f_A$, $f_A\in D(A)$,
$z\in \mathbb{C}\backslash \mathbb{R}$.
By~(\ref{f2_11}) we get
$$ Bx = BAf_A - zBf_A = ABf_A - zBf_A = (A-zE_H)Bf_A\in \Delta_A(z). $$
In particular, we have $BH_1\subseteq H_1$, $BH_3\subseteq H_3$.
Notice that $B_0^{-1}A_0 = A_0 B_0^{-1}$. It is a straightforward calculation to check that
\begin{equation}
\label{f2_20}
AB^{-1} x = B^{-1}A x,\qquad x\in D(A).
\end{equation}
Repeating the above argument with $B^{-1}$ instead of $B$ we get
$B^{-1}H_1\subseteq H_1$, $B^{-1}H_3\subseteq H_3$, and therefore
$H_1\subseteq BH_1$, $H_3\subseteq BH_3$. Consequently, the operator $B$ reduces
subspaces $H_1$ and $H_3$. It follows directly that $B$ reduces $H_2$ and $H_4$, as well.

\noindent
Since
$$ (A-iE_H) Bx = B(A-iE_H)x,\qquad x\in D(A), $$
for arbitrary $y\in H_1$, $y = (A-iE_H)x_A$, $x_A\in D(A)$, we have
$$ (A-iE_H) B (A-iE_H)^{-1} y = B y; $$
$$ B (A-iE_H)^{-1} y = (A-iE_H)^{-1} B y,\qquad y\in H_1, $$
and~(\ref{f2_19}) follows.
$\Box$

Our aim here is to construct a unitary operator $U$ in $H$, $U\supset V_A$, which commutes with $B$.
Choose an arbitrary $x\in H$, $x= x_{H_1} + x_{H_2}$. For an operator $U$ of the required type
by~Proposition~\ref{p2_1} we could write:
$$ BU x = BV_Ax_{H_1} + BU x_{H_2} = V_ABx_{H_1} + BU x_{H_2}, $$
$$ UB x = UB x_{H_1} + UB x_{H_2} = V_ABx_{H_1} + UB x_{H_2}. $$
So, it is enough to find an isometric operator $U_{2,4}$ which maps $H_2$ onto $H_4$, and
commutes with $B$:
\begin{equation}
\label{f2_21}
B U_{2,4} x = U_{2,4}B x,\qquad x\in H_2.
\end{equation}
Moreover, all operators $U$ of the required type have the following form:
\begin{equation}
\label{f2_22}
U = V_A \oplus U_{2,4},
\end{equation}
where $U_{2,4}$ is an isometric operator which maps $H_2$ onto $H_4$, and
commutes with $B$.

\noindent
We shall denote the operator $B$ restricted to $H_i$ by $B_{H_i}$, $1\leq i\leq 4$.
Notice that
\begin{equation}
\label{f2_23}
A^* J x= JA^* x,\qquad x\in D(A^*).
\end{equation}
Indeed, for arbitrary $f_A\in D(A)$ and $g_{A^*}\in D(A^*)$ we can write
$$ \overline{ (Af_A,Jg_{A^*})  } = (JAf_A, g_{A^*}) = (AJf_A, g_{A^*}) = (Jf_A, A^*g_{A^*}) $$
$$ = \overline{ (f_A,JA^*g_{A^*})  }, $$
and~(\ref{f2_23}) follows.

\noindent
Choose an arbitrary $x\in H_2$. We have
$$ A^* x = -i x, $$
and therefore
$$ A^* Jx = JA^* x = ix. $$
Thus, we have
$$ JH_2 \subseteq H_4. $$
In a similar manner we get
$$ JH_4 \subseteq H_2, $$
and therefore
\begin{equation}
\label{f2_24}
JH_2 = H_4,\quad JH_4 = H_2.
\end{equation}
By the Godi\v{c}-Lucenko Theorem (\cite{cit_2000_GL},\cite[Theorem 1]{cit_3000_GP}) we have a
representation:
\begin{equation}
\label{f2_25}
B_{H_2} = KL,
\end{equation}
where $K$ and $L$ are some conjugations in $H_2$.
We set
\begin{equation}
\label{f2_26}
U_{2,4} := JK.
\end{equation}
From~(\ref{f2_24}) it follows that $U_{2,4}$ maps isometrically $H_2$ onto $H_4$.
Notice that
\begin{equation}
\label{f2_27}
U_{2,4}^{-1} := KJ.
\end{equation}
Using relation~(\ref{f2_15}) we get
$$ U_{2,4} B_{H_2} U_{2,4}^{-1} x = JK KL KJ x = J LK J x = J B_{H_2}^{-1} J x $$
$$ = JB^{-1}J x = B x = B_{H_4} x,\qquad x\in H_4. $$
Therefore relation~(\ref{f2_21}) is true.

We define an operator $U$ by~(\ref{f2_22}) and define
\begin{equation}
\label{f2_28}
A_U := i(U+E_H)(U-E_H)^{-1} = iE_H + 2i(U-E_H)^{-1}.
\end{equation}
The inverse Cayley transformation $A_U$ is correctly defined since $1$ is not in the point spectrum of $U$.
Indeed, $V_A$ is the Cayley transformation of a symmetric operator while eigen subspaces $H_2$ and
$H_4$ have the zero intersection.
Let
\begin{equation}
\label{f2_29}
A_U = \int_\mathbb{R} s dE(s),\quad B = \int_{ [-\pi,\pi) } e^{i\varphi} dF(\varphi),
\end{equation}
where $E(s)$ and $F(\varphi)$ are the spectral measures of $A_U$ and $B$, respectively. These measures are
defined on $\mathfrak{B}(\mathbb{R})$ and $\mathfrak{B}([-\pi,\pi))$, respectively (\cite{cit_9000_BS}).
Since $U$ and $B$ commute, we get that $E(s)$ and $F(\varphi)$ commute, as well.
By induction argument we have
$$ x_{m,n} = A^m x_{0,n},\qquad m\in \mathbb{Z}_+,\ n\in \mathbb{Z}, $$
and
$$ x_{0,n} = B^n x_{0,0},\qquad n\in \mathbb{Z}. $$
Therefore we have
\begin{equation}
\label{f2_30}
x_{m,n} = A^m B^n x_{0,0},\qquad m\in \mathbb{Z}_+,\ n\in \mathbb{Z}.
\end{equation}
We can write
$$ x_{m,n} = \int_\mathbb{R} s^m dE(s) \int_{ [-\pi,\pi) } e^{in\varphi} dF(\varphi) x_{0,0} =
\int_\Pi s^m e^{in\varphi} d(E\times F) x_{0,0}, $$
where $E\times F$ is the product spectral measure on $\mathfrak{B}(\Pi)$.
Then
\begin{equation}
\label{f2_31}
s_{m,n} = (x_{m,n},x_{0,0})_H = \int_\Pi s^m e^{in\varphi} d((E\times F) x_{0,0}, x_{0,0})_H,\quad
(m,n)\in T.
\end{equation}
The measure $\mu := ((E\times F) x_{0,0}, x_{0,0})_H$ is a non-negative Borel measure on $\Pi$ and
relation~(\ref{f2_31}) shows that $\mu$ is a solution of the Devinatz moment problem.

Thus, we obtained a new proof of the following criterion.
\begin{thm}
\label{t2_2}
Let a Devinatz moment problem~(\ref{f1_1}) be given.
This problem has a solution if an only if conditions~(\ref{f2_1}) hold
for arbitrary complex numbers $\alpha_{m,n}$ such that all but finite numbers are zeros.
\end{thm}
{\bf Remark. } The original proof of Devinatz used the theory of reproducing kernels Hilbert spaces (RKHS).
In particular, he used properties of RKHS corresponding to the product of two positive matrices and
an inner structure of a RKHS corresponding to the moment problem. We used an abstract approach with
the Godi\v{c}-Lucenko Theorem and basic facts from the standard operator theory.

\section{Canonical solutions. A set of all solutions.}
Let a moment problem~(\ref{f1_1}) be given. Construct a Hilbert space $H$ and operators
$A,B,J$ as in the previous Section.
Let $\widetilde A\supseteq A$ be a self-adjoint extension of $A$ in a Hilbert space
$\widetilde H\supseteq H$. Let $R_z(\widetilde A)$, $z\in \mathbb{C}\backslash \mathbb{R}$, be the  resolvent
function of $\widetilde A$, and $E_{\widetilde A}$ be its spectral measure.
Recall that the function
\begin{equation}
\label{f3_1}
\mathbf{R}_z(A) := P^{\widetilde H}_H R_z(\widetilde A),\qquad z\in \mathbb{C}\backslash \mathbb{R},
\end{equation}
is said to be a generalized resolvent of $A$. The function
\begin{equation}
\label{f3_2}
\mathbf{E}_A (\delta) := P^{\widetilde H}_H E_{\widetilde A} (\delta),\qquad \delta\in \mathfrak{B}(\mathbb{R}),
\end{equation}
is said to be a spectral measure  of $A$.
There exists a one-to-one correspondence between generalized resolvents and spectral measures established by
the following relation~\cite{cit_6000_AG}:
\begin{equation}
\label{f3_3}
(\mathbf{R}_z(A) x,y)_H = \int_{\mathbb{R}} \frac{1}{t-z} d(\mathbf{E}_A x,y)_H,\qquad x,y\in H.
\end{equation}
We shall reduce the Devinatz moment problem to a problem of finding of generalized resolvents of a
certain class.
\begin{thm}
\label{t3_1}
Let a Devinatz moment problem~(\ref{f1_1}) be given and conditions~(\ref{f2_1}) hold.
Consider a Hilbert space $H$ and a sequence
$\{ x_{m,n} \}_{m\in \mathbb{Z}_+, n\in \mathbb{Z}}$, $x_{m,n}\in H$, such that relation~(\ref{f2_6})
holds where $K$ is defined by~(\ref{f2_2}).
Consider operators $A_0$,$B_0$ defined by~(\ref{f2_7}),(\ref{f2_8}) on
$L = \mathop{\rm Lin}\nolimits\{ x_{m,n} \}_{(m,n)\in T}$. Let $A=\overline{A_0}$, $B=\overline{B_0}$.
Let $\mu$ be an arbitrary solution of the moment problem. Then  it has the following form:
\begin{equation}
\label{f3_4}
\mu (\delta)= ((\mathbf{E}\times F)(\delta) x_{0,0}, x_{0,0})_H,\qquad \delta\in \mathfrak{B}(\mathbb{R}),
\end{equation}
where $F$ is the spectral measure of $B$, $\mathbf{E}$ is a spectral measure of $A$ which commutes with
$F$. By $((\mathbf{E}\times F)(\delta) x_{0,0}, x_{0,0})_H$ we mean the non-negative Borel measure on
$\mathbb{R}$ which is obtained by the Lebesgue continuation procedure from the following
non-negative measure on rectangules
\begin{equation}
\label{f3_5}
((\mathbf{E}\times F)(I_x\times I_\varphi)) x_{0,0}, x_{0,0})_H :=
( \mathbf{E}(I_x) F(I_\varphi)) x_{0,0}, x_{0,0})_H,
\end{equation}
where $I_x\subset \mathbb{R}$, $I_\varphi\subseteq [-\pi,\pi)$ are arbitrary intervals.

\noindent
On the other hand, for an arbitrary spectral measure $\mathbf{E}$ of $A$ which commutes with the
spectral measure $F$ of $B$, by relation~(\ref{f3_4}) it corresponds a solution of the moment
problem~(\ref{f1_1}).

\noindent
Moreover, the correspondence between the spectral measures of $A$ which commute with the spectral meeasure of
$B$ and solutions of the Devinatz moment problem is bijective.
\end{thm}
{\bf Remark. } The measure in~(\ref{f3_5}) is non-negative. Indeed,
for arbitrary intervals $I_x\subset \mathbb{R}$, $I_\varphi\subseteq [-\pi,\pi)$, we can write
$$ \left( \mathbf{E}(I_x) F(I_\varphi) x_{0,0}, x_{0,0} \right)_H =
\left( F(I_\varphi) \mathbf{E}(I_x) F(I_\varphi) x_{0,0}, x_{0,0} \right)_H $$
$$ = \left( \mathbf{E}(I_x) F(I_\varphi) x_{0,0}, F(I_\varphi)
x_{0,0} \right)_H = \left( \widehat E(I_x) F(I_\varphi) x_{0,0}, \widehat E(I_x) F(I_\varphi)
x_{0,0} \right)_{\widehat H} \geq 0, $$
where $\widehat E$ is the spectral function of a self-adjoint extension $\widehat A\supseteq A$ in
a Hilbert space $\widehat H\supseteq H$ such that $\mathbf{E} = P^{\widehat H}_H \widehat E$.
The measure in~(\ref{f3_5}) is additive. If $I_\varphi = I_{1,\varphi}\cup I_{2,\varphi}$,
$I_{1,\varphi}\cap I_{2,\varphi} = \emptyset$, then
$$ \left( \mathbf{E}(I_x) F(I_\varphi) x_{0,0}, x_{0,0} \right)_H =
\left( F( I_{1,\varphi}\cup I_{2,\varphi} )\mathbf{E}(I_x) x_{0,0}, x_{0,0} \right)_H $$
$$ = \left( F(I_{1,\varphi})\mathbf{E}(I_x) x_{0,0}, x_{0,0} \right)_H +
\left( F(I_{2,\varphi})\mathbf{E}(I_x) x_{0,0}, x_{0,0} \right)_H. $$
The case $I_x = I_{1,x}\cup I_{2,x}$ is analogous.
Moreover, repeating the standard arguments~\cite[Chapter 5, Theorem 2, p. 254-255]{cit_9500_KF} we conclude
that the measure in~(\ref{f3_5}) is $\sigma$-additive.
Thus, it posesses the (unique) Lebesgue continuation to a (finite) non-negative Borel measure
on $\Pi$.

{\bf Proof. }
Consider a Hilbert space $H$ and operators $A$,$B$ as in the statement of the Theorem.
Let $F$ be the spectral measure of $B$. Let $\mu$ be an arbitrary solution of the moment problem~(\ref{f1_1}).
Consider the space $L^2_\mu$ of complex functions on $\Pi$ which are square integrable with respect to the
measure $\mu$. The scalar product and the norm are given by
$$ (f,g)_\mu =
\int_\Pi f(x,\varphi) \overline{ g(x,\varphi) } d\mu,\quad
\|f\|_\mu = \left( (f,f)_\mu \right)^{ \frac{1}{2} },\quad f,g\in L^2_\mu. $$
Consider the following operators:
\begin{equation}
\label{f3_6}
A_\mu f(x,\varphi) = xf(x,\varphi),\qquad D(A_\mu) = \{ f\in L^2_\mu:\ xf(x,\varphi)\in L^2_\mu \},
\end{equation}
\begin{equation}
\label{f3_7}
B_\mu f(x,\varphi) = e^{i\varphi} f(x,\varphi),\qquad D(B_\mu) = L^2_\mu.
\end{equation}
The operator $A_\mu$ is self-adjoint and the operator $B_\mu$ is unitary. Moreover, these operators
commute and therefore the spectral measure $E_\mu$ of $A_\mu$ and the spectral measure $F_\mu$ of $B_\mu$
commute, as well.

\noindent
Let $p(x,\varphi)$ be a (power-trigonometric) polynomial of the form~(\ref{f1_1}) and
$q(x,\varphi)$ be a (power-trigonometric) polynomial of the form~(\ref{f1_1}) with
$\beta_{m,n}\in \mathbb{C}$ instead of $\alpha_{m,n}$.
Then
$$ (p,q)_\mu = \sum_{(m,n)\in T, (k,l)\in T} \alpha_{m,n}\overline{ \beta_{k,l} }
\int_\Pi x^{m+k} e^{i(n-l)\varphi} d\mu $$
$$ = \sum_{(m,n)\in T, (k,l)\in T} \alpha_{m,n}\overline{ \beta_{k,l} } s_{m+k,n-l}, $$
On the other hand, we can write
$$ \left(
\sum_{(m,n)\in T} \alpha_{m,n} x_{m,n}, \sum_{(k,l)\in T} \beta_{k,l} x_{k,l} \right)_H =
\sum_{(m,n)\in T, (k,l)\in T} \alpha_{m,n}\overline{ \beta_{k,l} }
(x_{m,n},x_{k,l})_H $$
$$ = \sum_{(m,n)\in T, (k,l)\in T} \alpha_{m,n}\overline{ \beta_{k,l} } K((m,n),(k,l))
= \sum_{(m,n)\in T, (k,l)\in T} \alpha_{m,n}\overline{ \beta_{k,l} } s_{m+k,n-l}. $$
Therefore
\begin{equation}
\label{f3_8}
(p,q)_\mu = \left(
\sum_{(m,n)\in T} \alpha_{m,n} x_{m,n}, \sum_{(k,l)\in T} \beta_{k,l} x_{k,l} \right)_H.
\end{equation}
Consider thr following operator:
\begin{equation}
\label{f3_9}
V[p] = \sum_{(m,n)\in T} \alpha_{m,n} x_{m,n},\quad p=\sum_{(m,n)\in T} \alpha_{m,n} x^m e^{in\varphi}.
\end{equation}
Here by $[p]$ we mean the class of equivalence in $L^2_\mu$ defined by $p$. If two different polynomials
$p$ and $q$ belong to the same class of equivalence then by~(\ref{f3_8}) we get
$$ 0 = \| p-q \|_\mu^2 = (p-q,p-q)_\mu = \left( \sum_{(m,n)\in T} (\alpha_{m,n}-\beta_{m,n}) x_{m,n},
\sum_{(k,l)\in T} (\alpha_{k,l}-\beta_{k,l}) x_{k,l} \right) $$
$$ = \left\| \sum_{(m,n)\in T} \alpha_{m,n} x_{m,n} - \sum_{(m,n)\in T} \beta_{m,n} x_{m,n} \right\|_\mu^2. $$
Thus, the definition of $V$ is correct. It is not hard to see that $V$ maps a set of polynomials $P^2_{0,\mu}$
in $L^2_\mu$ on $L$. By continuity we extend $V$ to the isometric transformation from
the closure of polynomials $P^2_\mu = \overline{P^2_{0,\mu}}$ onto $H$.

\noindent
Set $H_0 := L^2_\mu \ominus P^2_\mu$. Introduce the following operator:
\begin{equation}
\label{f3_10}
U := V \oplus E_{H_0},
\end{equation}
which maps isometrically $L^2_\mu$ onto $\widetilde H := H\oplus H_0$.
Set
\begin{equation}
\label{f3_11}
\widetilde A := UA_\mu U^{-1},\quad \widetilde B := UB_\mu U^{-1}.
\end{equation}
Notice that
$$ \widetilde A x_{m,n} = UA_\mu U^{-1} x_{m,n} = UA_\mu x^m e^{in\varphi} = Ux^{m+1} e^{in\varphi} = x_{m+1,n}, $$
$$ \widetilde B x_{m,n} = UB_\mu U^{-1} x_{m,n} = UB_\mu x^m e^{in\varphi} = Ux^{m} e^{i(n+1)\varphi} = x_{m,n+1}. $$
Therefore $\widetilde A\supseteq A$ and $\widetilde B\supseteq B$.
Let
\begin{equation}
\label{f3_12}
\widetilde A = \int_\mathbb{R} s d\widetilde E(s),\quad \widetilde B = \int_{ [-\pi,\pi) } e^{i\varphi}
d \widetilde F(\varphi),
\end{equation}
where $\widetilde E(s)$ and $\widetilde F(\varphi)$ are the spectral measures of $\widetilde A$ and
$\widetilde B$, respectively.
Repeating arguments after relation~(\ref{f2_29}) we obtain that
\begin{equation}
\label{f3_13}
x_{m,n} = \widetilde A^m \widetilde B^n x_{0,0},\qquad m\in \mathbb{Z}_+,\ n\in \mathbb{Z},
\end{equation}
\begin{equation}
\label{f3_14}
s_{m,n} = \int_\Pi s^m e^{in\varphi} d((\widetilde E\times \widetilde F) x_{0,0}, x_{0,0})_{\widetilde H},\quad
(m,n)\in T,
\end{equation}
where $(\widetilde E\times \widetilde F)$ is the product measure of $\widetilde E$ and $\widetilde F$.
Thus, the measure $\widetilde \mu := ((\widetilde E\times \widetilde F) x_{0,0}, x_{0,0})_{\widetilde H}$
is a solution of the Devinatz moment problem.

\noindent
Let $I_x\subset \mathbb{R}$, $I_\varphi\subseteq [-\pi,\pi)$ be arbitrary intervals.
Then
$$ \widetilde \mu (I_x \times I_\varphi) = ((\widetilde E\times \widetilde F) (I_x \times I_\varphi)
x_{0,0}, x_{0,0})_{\widetilde H} $$
$$ = ( \widetilde E(I_x) \widetilde F(I_\varphi) x_{0,0}, x_{0,0})_{\widetilde H} =
( P^{\widetilde H}_H \widetilde E(I_x) \widetilde F(I_\varphi) x_{0,0}, x_{0,0})_{\widetilde H} $$
$$ = ( \mathbf{E}(I_x) F(I_\varphi) x_{0,0}, x_{0,0})_{H}, $$
where $\mathbf{E}$ is the correponding spectral function of $A$ and $F$ is the spectral function of $B$.
Thus, the measure $\widetilde \mu$ has the form~(\ref{f3_4}) since the Lebesgue continuation
is unique.

\noindent
Let us show that $\widetilde \mu = \mu$.
Consider the following transformation:
\begin{equation}
\label{f3_15}
S:\ (x,\varphi) \in \Pi \mapsto \left( \mathop{\rm Arg }\nolimits \frac{x-i}{x+i}, \varphi  \right) \in \Pi_0,
\end{equation}
where $\Pi_0 = [-\pi,\pi) \times [-\pi,\pi)$ and $\mathop{\rm Arg }\nolimits e^{iy} = y\in [-\pi,\pi)$.
By virtue of $V$ we define  the following measures:
\begin{equation}
\label{f3_16}
\mu_0 (VG) := \mu (G),\quad \widetilde\mu_0 (VG) := \widetilde\mu (G),\qquad G\in \mathfrak{B}(\Pi),
\end{equation}
It is not hard to see that $\mu_0$ and $\widetilde\mu_0$ are non-negative measures on
$\mathfrak{B}(\Pi_0)$.
Then
\begin{equation}
\label{f3_17}
\int_\Pi \left( \frac{x-i}{x+i} \right)^m e^{in\varphi} d\mu =
\int_{\Pi_0} e^{im\psi} e^{in\varphi} d\mu_0,
\end{equation}
\begin{equation}
\label{f3_18}
\int_\Pi \left( \frac{x-i}{x+i} \right)^m e^{in\varphi} d\widetilde\mu =
\int_{\Pi_0} e^{im\psi} e^{in\varphi} d\widetilde\mu_0,\qquad m,n\in \mathbb{Z};
\end{equation}
and
$$  \int_\Pi \left( \frac{x-i}{x+i} \right)^m e^{in\varphi} d\widetilde\mu =
\int_\Pi \left( \frac{x-i}{x+i} \right)^m e^{in\varphi}
d((\widetilde E\times \widetilde F) x_{0,0}, x_{0,0})_{\widetilde H} $$
$$ = \left( \int_\Pi \left( \frac{x-i}{x+i} \right)^m e^{in\varphi}
d(\widetilde E\times \widetilde F) x_{0,0}, x_{0,0} \right)_{\widetilde H}  $$
$$ = \left( \int_\mathbb{R} \left( \frac{x-i}{x+i} \right)^m d\widetilde E
\int_{[-\pi,\pi)} e^{in\varphi} d\widetilde F x_{0,0}, x_{0,0} \right)_{\widetilde H}  $$
$$ = \left( \left( (\widetilde A - iE_{\widetilde H})(\widetilde A + iE_{\widetilde H})^{-1} \right)^m
\widetilde B^n  x_{0,0}, x_{0,0} \right)_{\widetilde H}  $$
$$ = \left( U^{-1}\left( (\widetilde A - iE_{\widetilde H})(\widetilde A + iE_{\widetilde H})^{-1} \right)^m
\widetilde B^n  U 1, U 1 \right)_\mu  $$
$$ = \left( \left( (A_\mu - iE_{L^2_\mu})(A_\mu + iE_{L^2_\mu})^{-1} \right)^m
B_\mu^n  1, 1 \right)_\mu  $$
\begin{equation}
\label{f3_19}
 = \int_\Pi \left( \frac{x-i}{x+i} \right)^m e^{in\varphi} d\mu,\qquad m,n\in \mathbb{Z}.
\end{equation}
By virtue of relations~(\ref{f3_17}),(\ref{f3_18}) and~(\ref{f3_19}) we get
\begin{equation}
\label{f3_20}
\int_{\Pi_0} e^{im\psi} e^{in\varphi} d\mu_0 =
\int_{\Pi_0} e^{im\psi} e^{in\varphi} d\widetilde\mu_0,\qquad m,n\in \mathbb{Z}.
\end{equation}
By the Weierstrass theorem we can approximate any continuous function by exponentials and therefore
\begin{equation}
\label{f3_21}
\int_{\Pi_0} f(\psi) g(\varphi) d\mu_0 =
\int_{\Pi_0} f(\psi) g(\varphi) d\widetilde\mu_0,
\end{equation}
for arbitrary continuous functions on $\Pi_0$. In particular, we have
\begin{equation}
\label{f3_22}
\int_{\Pi_0} \psi^n \varphi^m d\mu_0 =
\int_{\Pi_0} \psi^n \varphi^m d\widetilde\mu_0,\qquad n,m\in \mathbb{Z}_+.
\end{equation}
However, the two-dimensional Hausdorff moment problem is determinate (\cite{cit_10000_ST}) and therefore we get
$\mu_0 = \widetilde\mu_0$ and $\mu=\mu_0$.
Thus, we have proved that an arbitrary solution $\mu$ of the Devinatz moment problem can be represented
in the form~(\ref{f3_4}).

Let us check the second assertion of the Theorem.
For an arbitrary spectral measure $\mathbf{E}$ of $A$ which commutes with the
spectral measure $F$ of $B$, by relation~(\ref{f3_4}) we define a non-negative Borel measure $\mu$
on $\Pi$. Let us show that the measure $\mu$ is a solution of the moment
problem~(\ref{f1_1}).

\noindent
Let $\widehat A$ be a self-adjoint extension of the operator $A$ in a Hilbert space
$\widehat H\supseteq H$, such that
$$ \mathbf{E} = P^{\widehat H}_H \widehat E, $$
where $\widehat E$ is the spectral measure of $\widehat A$.
By~(\ref{f2_30}) we get
$$ x_{m,n} = A^m B^n x_{0,0} = \widehat A^m B^n x_{0,0} = P^{\widehat H}_H \widehat A^m B^n x_{0,0} $$
$$ = P^{\widehat H}_H \left( \lim_{a\to +\infty} \int_{[-a,a)} x^m d\widehat E \right)
\int_{[-\pi,\pi)} e^{in\varphi} dF x_{0,0}
= \left( \lim_{a\to +\infty} \int_{[-a,a)} x^m d\mathbf{E} \right) $$
$$ * \int_{[-\pi,\pi)} e^{in\varphi} dF x_{0,0}
= \left( \lim_{a\to +\infty} \left( \int_{[-a,a)} x^m d\mathbf{E} \int_{[-\pi,\pi)} e^{in\varphi} dF \right)
\right) x_{0,0}, $$
\begin{equation}
\label{f3_23}
\qquad m\in \mathbb{Z}_+,\ n\in \mathbb{Z},
\end{equation}
where the limits are understood in the weak operator topology.
Then we choose arbitrary points
$$  -a = x_0 < x_1 < ... < x_{N}=a; $$
\begin{equation}
\label{f3_24}
\max_{1\leq i\leq N}|x_{i}-x_{i-1}| =: d,\quad  N\in \mathbb{N};
\end{equation}
$$  -\pi = \varphi_0 < \varphi_1 < ... < \varphi_{M}=\pi; $$
\begin{equation}
\label{f3_25}
 \max_{1\leq j\leq M}|\varphi_{j}-\varphi_{j-1}| =: r;\quad M\in \mathbb{N}.
\end{equation}
Set
$$ C_a := \int_{[-a,a)} x^m d\mathbf{E} \int_{[-\pi,\pi)} e^{in\varphi} dF =
\lim_{d\rightarrow 0} \sum_{i=1}^N x_{i-1}^m \mathbf{E}([x_{i-1},x_i)) $$
$$ * \lim_{r\rightarrow 0} \sum_{j=1}^M e^{in\varphi_{j-1}} F([\varphi_{j-1},\varphi_j)), $$
where the integral sums converge in the strong operator topology. Then
$$ C_a = \lim_{d\rightarrow 0} \lim_{r\rightarrow 0} \sum_{i=1}^N x_{i-1}^m \mathbf{E}([x_{i-1},x_i))
\sum_{j=1}^M e^{in\varphi_{j-1}} F([\varphi_{j-1},\varphi_j)) $$
$$ = \lim_{d\rightarrow 0} \lim_{r\rightarrow 0}
\sum_{i=1}^N \sum_{j=1}^M
x_{i-1}^m e^{in\varphi_{j-1}}
\mathbf{E}([x_{i-1},x_i)) F([\varphi_{j-1},\varphi_j)), $$
where the limits are understood in the strong operator topology. Then
$$ (C_a x_{0,0}, x_{0,0})_H =
\left( \lim_{d\rightarrow 0} \lim_{r\rightarrow 0}
\sum_{i=1}^N \sum_{j=1}^M
x_{i-1}^m e^{in\varphi_{j-1}}
\mathbf{E}([x_{i-1},x_i)) F([\varphi_{j-1},\varphi_j)) x_{0,0}, x_{0,0} \right)_H $$
$$ = \lim_{d\rightarrow 0} \lim_{r\rightarrow 0}
\sum_{i=1}^N \sum_{j=1}^M
x_{i-1}^m e^{in\varphi_{j-1}}
\left( \mathbf{E}([x_{i-1},x_i)) F([\varphi_{j-1},\varphi_j)) x_{0,0}, x_{0,0} \right)_H $$
$$ = \lim_{d\rightarrow 0} \lim_{r\rightarrow 0}
\sum_{i=1}^N \sum_{j=1}^M
x_{i-1}^m e^{in\varphi_{j-1}}
\left( (\mathbf{E}\times F) ( [x_{i-1},x_i)\times [\varphi_{j-1},\varphi_j) ) x_{0,0}, x_{0,0} \right)_H $$
$$ = \lim_{d\rightarrow 0} \lim_{r\rightarrow 0}
\sum_{i=1}^N \sum_{j=1}^M
x_{i-1}^m e^{in\varphi_{j-1}}
\left( \mu ( [x_{i-1},x_i)\times [\varphi_{j-1},\varphi_j) ) x_{0,0}, x_{0,0} \right)_H. $$
Therefore
$$ (C_a x_{0,0}, x_{0,0})_H =
\lim_{d\rightarrow 0} \lim_{r\rightarrow 0}
\int_{[-a,a)\times[-\pi,\pi)} f_{d,r} (x,\varphi) d\mu, $$
where $f_{d,r}$ is equal to $x_{i-1}^m e^{in\varphi_{j-1}}$ on the rectangular
$[x_{i-1},x_i) \times [\varphi_{j-1},\varphi_j)$, $1\leq i\leq N$, $1\leq j\leq M$.

\noindent
If $r\rightarrow 0$, then the simple function
$f_{d,r}$ converges uniformly to the function $f_d$ which is equal to
$x_{i-1}^m e^{in\varphi}$ on the rectangular
$[x_{i-1},x_i) \times [\varphi_{j-1},\varphi_j)$, $1\leq i\leq N$, $1\leq j\leq M$.
Then
$$ (C_a x_{0,0}, x_{0,0})_H =
\lim_{d\rightarrow 0}
\int_{[-a,a)\times[-\pi,\pi)} f_{d} (x,\varphi) d\mu. $$
If $d\rightarrow 0$, then the function
$f_{d}$ converges uniformly to the function $x^m e^{in\varphi}$. Since
$|f_d|\leq A^m$, by the Lebesgue theorem we get
\begin{equation}
\label{f3_26}
(C_a x_{0,0}, x_{0,0})_H =
\int_{[-a,a)\times[-\pi,\pi)} x^m e^{in\varphi} d\mu.
\end{equation}
By virtue of relations~(\ref{f3_23}) and~(\ref{f3_26}) we get
$$ s_{m,n} = (x_{m,n},x_{0,0})_H = \lim_{a\to +\infty} (C_a x_{0,0},x_{0,0})_H $$
\begin{equation}
\label{f3_27}
= \lim_{a\to+\infty} \int_{[-a,a)\times[-\pi,\pi)} x^m e^{in\varphi} d\mu =
\int_\Pi x^m e^{in\varphi} d\mu.
\end{equation}
Thus, the measure $\mu$ is a solution of the Devinatz moment problem.

Let us prove the last assertion of the Theorem. Suppose to the contrary that two different
spectral measures $\mathbf{E}_1$ and $\mathbf{E}_1$ of $A$ commute with the spectral measure $F$ of
$B$ and produce by relation~(\ref{f3_4}) the same solution $\mu$ of the Devinatz moment problem.
Choose an arbitrary $z\in \mathbb{C}\backslash \mathbb{R}$. Then
$$ \int_\Pi \frac{x^m}{x-z} e^{in\varphi} d\mu =
\int_\Pi \frac{x^m}{x-z} e^{in\varphi} ((\mathbf{E}_k\times F)(\delta) x_{0,0}, x_{0,0})_H $$
\begin{equation}
\label{f3_28}
= \lim_{a\to +\infty}
\int_{[-a,a)\times [-\pi,\pi)}
\frac{x^m}{x-z} e^{in\varphi} d((\mathbf{E}_k\times F)(\delta) x_{0,0}, x_{0,0})_H,\quad k=1,2.
\end{equation}
Consider arbitrary partitions of the type~(\ref{f3_24}),(\ref{f3_25}). Then
$$ D_a := \int_{[-a,a)\times [-\pi,\pi)}
\frac{x^m}{x-z} e^{in\varphi} d((\mathbf{E}_k\times F)(\delta) x_{0,0}, x_{0,0})_H $$
$$ = \lim_{d\to 0} \lim_{r\to 0}
\int_{[-a,a)\times [-\pi,\pi)} g_{z;d,r}(x,\varphi)
d((\mathbf{E}_k\times F)(\delta) x_{0,0}, x_{0,0})_H. $$
Here the function $g_{z;d,r}(x,\varphi)$ is equal to
$\frac{x_{i-1}^m}{x_{i-1}-z} e^{in\varphi_{j-1}}$
on the rectangular
$[x_{i-1},x_i) \times [\varphi_{j-1},\varphi_j)$, $1\leq i\leq N$, $1\leq j\leq M$.
Then
$$ D_a = \lim_{d\to 0} \lim_{r\to 0}
\sum_{i=1}^N \sum_{j=1}^M
\frac{ x_{i-1}^m }{ x_{i-1}-z } e^{in\varphi_{j-1}}
\left( \mathbf{E}_k ([x_{i-1},x_i)) F([\varphi_{j-1},\varphi_j)) x_{0,0}, x_{0,0} \right)_H $$
$$ =
\lim_{d\to 0} \lim_{r\to 0}
\left( \sum_{i=1}^N
\frac{ x_{i-1}^m }{ x_{i-1}-z }  \mathbf{E}_k ([x_{i-1},x_i))
\sum_{j=1}^M
e^{in\varphi_{j-1}} F([\varphi_{j-1},\varphi_j)) x_{0,0}, x_{0,0} \right)_H $$
$$ =
\left( \int_{[-a,a)}
\frac{ x^m }{ x-z }  d\mathbf{E}_k \int_{[-\pi,\pi)}
e^{in\varphi} dF x_{0,0}, x_{0,0} \right)_H. $$
Let $n = n_1+n_2$, $n_1,n_2\in \mathbb{Z}$. Then we can write:
$$ D_a = \left( B^{n_1} \int_{[-a,a)}
\frac{ x^m }{ x-z }  d\mathbf{E}_k B^{n_2} x_{0,0}, x_{0,0} \right)_H $$
$$ =
\left( \int_{[-a,a)}
\frac{ x^m }{ x-z }  d\mathbf{E}_k x_{0,n_2}, x_{0,-n_1} \right)_H. $$
By~(\ref{f3_28}) we get
$$ \int_\Pi \frac{x^m}{x-z} e^{in\varphi} d\mu =
\lim_{a\to +\infty} D_a =
\lim_{a\to +\infty}\left( \int_{[-a,a)}
\frac{ x^m }{ x-z }  d \widehat{E}_k x_{0,n_2}, x_{0,-n_1} \right)_{\widehat H_k} $$
$$ = \left( \int_\mathbb{R}
\frac{ x^m }{ x-z }  d\widehat{E}_k x_{0,n_2}, x_{0,-n_1} \right)_{\widehat H_k}
= \left( \widehat{A}^{m_2} R_z(\widehat{A}_k) \widehat{A}^{m_1} x_{0,n_2}, x_{0,-n_1} \right)_{\widehat H_k}
$$
\begin{equation}
\label{f3_29}
= \left( R_z(\widehat{A}_k) x_{m_1,n_2}, x_{m_2,-n_1} \right)_H,
\end{equation}
where $m_1,m_2\in \mathbb{Z}_+:\ m_1+m_2 = m$,
and $\widehat A_k$ is a self-adjoint extension of $A$ in a Hilbert space $\widehat H_k\supseteq H$ such that
its spectral measure $\widehat E_k$ generates $\mathbf{E}_k$: $\mathbf{E}_k = P^{\widehat H_k}_H \widehat E_k$;
$k=1,2$.

\noindent
Relation~(\ref{f3_29}) shows that the generalized resolvents corresponding to $\mathbf{E}_k$, $k=1,2$, coincide.
That means that the spectral measures $\mathbf{E}_1$ and $\mathbf{E}_2$ coincide. We obtained a contradiction.
This completes the proof.
$\Box$

\begin{dfn}
\label{d3_1}
A solution $\mu$ of the Devinatz moment problem~(\ref{f1_1}) we shall call {\bf canonical}
if it is generated by relation~(\ref{f3_4}) where $\mathbf{E}$ is an {\bf orthogonal}
spectral measure of $A$ which commutes with the spectral measure of $B$. Orthogonal spectral measures
are those measures which are the spectral measures of self-adjoint extensions of $A$ inside $H$.
\end{dfn}
Let a moment problem~(\ref{f1_1}) be given and conditions~(\ref{f2_1}) hold.
Let us describe canonical solutions of the Devinatz moment problem.
In the proof of Theorem~\ref{t2_2} we have constructed one canonical solution, see relation~(\ref{f2_31}).
Let $\mu$ be an arbitrary canonical solution and $\mathbf{E}$ be the corresponding orthogonal spectral
measure of $A$. Let $\widetilde A$ be the self-adjoint operator in $H$ which corresponds to $\mathbf{E}$.
Consider the Cayley transformation of $\widetilde A$:
\begin{equation}
\label{f3_30}
U_{\widetilde A} = (\widetilde A + iE_H)(\widetilde A - iE_H)^{-1} \supseteq V_A,
\end{equation}
where $V_A$ is defined by~(\ref{f2_16}).
Since $\mathbf{E}$ commutes with the spectral measure $F$ of $B$, then $U_{\widetilde A}$ commutes
with $B$.
By relation~(\ref{f2_22}) the operator $U_{\widetilde A}$ have the following form:
\begin{equation}
\label{f3_31}
U_{\widetilde A} = V_A \oplus \widetilde U_{2,4},
\end{equation}
where $\widetilde U_{2,4}$ is an isometric operator which maps $H_2$ onto $H_4$, and
commutes with $B$.
Let the operator $U_{2,4}$ be defined by~(\ref{f2_26}). Then the following operator
\begin{equation}
\label{f3_32}
U_2 = U_{2,4}^{-1} \widetilde U_{2,4},
\end{equation}
is a unitary operator in $H_2$ which commutes with $B_{H_2}$.

Denote by $\mathbf{S}(B;H_2)$ a set of all unitary operators in $H_2$ which commute with $B_{H_2}$.
Choose an arbitrary operator $\widehat U_2\in \mathbf{S}(B;H_2)$. Define
$\widehat U_{2,4}$ by the following relation:
\begin{equation}
\label{f3_33}
\widehat U_{2,4} = U_{2,4} \widehat U_2.
\end{equation}
Notice that $\widehat U_{2,4}$ commutes with $B_{H_2}$.
Then we define a unitary operator $U = V_A \oplus \widehat U_{2,4}$ and its Cayley transformation
$\widehat A$ which commute with the operator $B$.
Repeating arguments before~(\ref{f2_31}) we get a canonical solution of the Devinatz moment problem.

\noindent
Thus, all canonical solutions of the Devinatz moment problem are generated by operators
$\widehat U_2\in \mathbf{S}(B;H_2)$. Notice that different operators $U',U''\in \mathbf{S}(B;H_2)$ produce different
orthogonal spectral measures $\mathbf{E}',\mathbf{E}$. By Theorem~\ref{t3_1},
these spectral measures produce different solutions of the moment problem.

Recall some definitions from~\cite{cit_9000_BS}.
A pair $(Y,\mathfrak{A})$, where $Y$ is an arbitrary set and $\mathfrak{A}$ is a fixed
$\sigma$-algebra of subsets of $A$ is said to be a {\it measurable space}.
A triple $(Y,\mathfrak{A},\mu)$, where $(Y,\mathfrak{A})$ is a measurable space and
$\mu$ is a measure on $\mathfrak{A}$ is said to be a {\it space with a measure}.

Let $(Y,\mathfrak{A})$ be  a measurable space, $\mathbf{H}$ be a Hilbert space and
$\mathcal{P}=\mathcal{P}(\mathbf{H})$ be a set of all orthogonal projectors in $\mathbf{H}$.
A countably additive mapping $E:\ \mathfrak{A}\rightarrow \mathcal{P}$, $E(Y) = E_{\mathbf{H}}$,
is said to be a {\it spectral measure} in $\mathbf{H}$.
A set $(Y,\mathfrak{A},H,E)$ is said to be a {\it space with a spectral measure}.
By $S(Y,E)$ one means a set of all $E$-measurable $E$-a.e. finite complex-valued functions on $Y$.

Let $(Y,\mathfrak{A},\mu)$ be a separable space with a $\sigma$-finite measure and
to $\mu$-everyone $y\in Y$ it corresponds a Hilbert space $G(y)$. A function
$N(y) = \dim G(y)$ is called the {\it dimension function}.
It is supposed to be $\mu$-measurable. Let $\Omega$ be a set of vector-valued functions $g(y)$ with
values in $G(y)$ which are defined $\mu$-everywhere and are measurable with respect to some
base of measurability. A set of (classes of equivalence) of such functions with
the finite norm
\begin{equation}
\label{f3_34}
\| g \|^2_{\mathcal{H}} = \int |g(y)|^2_{G(y)} d\mu(y) <\infty
\end{equation}
form a Hilbert space $\mathcal{H}$ with the scalar product given by
\begin{equation}
\label{f3_35}
( g_1,g_2 )_{\mathcal{H}} = \int (g_1,g_2)_{G(y)} d\mu(y).
\end{equation}
The space $\mathcal{H}= \mathcal{H}_{\mu,N} = \int_Y \oplus G(y) d\mu(y)$
is said to be a {\it direct integral of Hilbert spaces}.
Consider the following operator
\begin{equation}
\label{f3_36}
\mathbf{X}(\delta) g = \chi_\delta g,\qquad g\in \mathcal{H},\ \delta\in \mathfrak{A},
\end{equation}
where $\chi_\delta$ is the characteristic function of the set $\delta$.
The operator $\mathbf{X}$ is a spectral measure in $\mathcal{H}$.

Let $t(y)$ be a measurable operator-valued function with values in $\mathbf{B}(G(y))$ which is
$\mu$-a.e. defined and $\mu-\sup \|t(y)\|_{G(y)} < \infty$. The operator
\begin{equation}
\label{f3_37}
T:\ g(y) \mapsto t(y)g(y),
\end{equation}
is said to be {\it decomposable}. It is a bounded operator in $\mathcal{H}$ which commutes with
$\mathbf{X}(\delta)$, $\forall\delta\in \mathfrak{A}$.
Moreover, every bounded operator in $\mathcal{H}$ which commutes with
$\mathbf{X}(\delta)$, $\forall\delta\in \mathfrak{A}$, is decomposable~\cite{cit_9000_BS}.
In the case $t(y) = \varphi(y)E_{G(y)}$, where $\varphi\in S(Y,\mu)$, we set $T =: Q_\varphi$.
The decomposable operator is unitary if and only if $\mu$-a.e. the operator $t(y)$ is unitary.

Return to the study of canonical solutions. Consider the spectral measure
$F_2$ of the operator $B_{H_2}$ in $H_2$.
There exists an element $h\in H_2$ of the maximal type, i.e. the non-negative Borel measure
\begin{equation}
\label{f3_38}
\mu(\delta) := (F_2(\delta)h,h),\qquad \delta\in \mathfrak{B}([-\pi,\pi)),
\end{equation}
has the maximal type between all such measures (generated by other elements of $H_2$). This type
is said to be the {\it spectral type} of the measure $F_2$.
Let $N_2$ be the multiplicity function of the measure $F_2$. Then there exists a unitary transformation $W$
of the space $H_2$ on $\mathcal{H}=\mathcal{H}_{\mu,N_2}$ such that
\begin{equation}
\label{f3_39}
W B_{H_2} W^{-1} = Q_{e^{iy}},\qquad W F_2(\delta) W^{-1} = \mathbf{X}(\delta).
\end{equation}
Notice that $\widehat U_2\in \mathbf{S}(B;H_2)$ if and only if
the operator
\begin{equation}
\label{f3_40}
V_2 := W \widehat U_2 W^{-1},
\end{equation}
is unitary and commutes with $\mathbf{X}(\delta)$, $\forall\delta\in \mathfrak{[-\pi,\pi)}$.
The latter is equivalent to the condition that $V_2$ is decomposable and the values of the corresponding
operator-valued function $t(y)$ are $\mu$-a.e. unitary operators.
A set of all decomposable operators in $\mathcal{H}$ such that the values of the corresponding
operator-valued function $t(y)$ are $\mu$-a.e. unitary operators we denote by $\mathbf{D}(B;H_2)$.

\begin{thm}
\label{t3_2}
Let a Devinatz moment problem~(\ref{f1_1}) be given. In conditions of Theorem~\ref{t3_1} all
canonical solutions of the moment problem have the  form~(\ref{f3_4}) where the spectral
measures $\mathbf{E}$ of the operator $A$ are constructed by operators from $\mathbf{D}(B;H_2)$.
Namely, for an arbitrary $V_2\in \mathbf{D}(B;H_2)$ we set $U_2 = W^{-1} V_2 W$,
$\widehat U_{2,4} = U_{2,4} \widehat U_2$, $U = V_A \oplus \widehat U_{2,4}$,
$\widehat A = i(U+E_H)(U-E_H)^{-1}$, and then $\mathbf{E}$ is the spectral measure of $\widehat A$.

\noindent
Moreover, the correspondence between $\mathbf{D}(B;H_2)$ and a set of all canonical solutions of
the Devinatz moment problem is bijective.
\end{thm}
{\bf Proof. } The proof follows directly from the previous considerations.
$\Box$

Consider a Devinatz moment problem~(\ref{f1_1}) and suppose that conditions~(\ref{f2_1}) hold.
Let us turn to a parameterization of all solutions of the moment problem.
We shall use Theorem~\ref{t3_1}. Consider relation~(\ref{f3_4}). The spectral measure $\mathbf{E}$
commutes with the operator $B$.
Choose an arbitrary $z\in \mathbb{C}\backslash \mathbb{R}$.
By virtue of relation~(\ref{f3_3}) we can write:
$$ (B\mathbf{R}_z(A) x,y)_H = (\mathbf{R}_z(A) x,B^*y)_H =
\int_{\mathbb{R}} \frac{1}{t-z} d(\mathbf{E}(t) x,B^*y)_H $$
\begin{equation}
\label{f3_41}
\int_{\mathbb{R}} \frac{1}{t-z} d(B\mathbf{E}(t) x,y)_H
= \int_{\mathbb{R}} \frac{1}{t-z} d(\mathbf{E}(t)B x,y)_H,\qquad x,y\in H;
\end{equation}
\begin{equation}
\label{f3_42}
(\mathbf{R}_z(A) Bx,y)_H = \int_{\mathbb{R}} \frac{1}{t-z} d(\mathbf{E}(t) Bx,y)_H,\qquad x,y\in H,
\end{equation}
where $\mathbf{R}_z(A)$ is  the generalized resolvent which corresponds to $\mathbf{E}$.
Therefore we get
\begin{equation}
\label{f3_43}
\mathbf{R}_z(A) B = B \mathbf{R}_z(A),\qquad z\in \mathbb{C}\backslash \mathbb{R}.
\end{equation}
On the other hand, if relation~(\ref{f3_43}) holds, then
\begin{equation}
\label{f3_44}
\int_{\mathbb{R}} \frac{1}{t-z} d(\mathbf{E} Bx,y)_H =
\int_{\mathbb{R}} \frac{1}{t-z} d(B\mathbf{E} x,y)_H,\quad x,y\in H,\ z\in \mathbb{C}\backslash \mathbb{R}.
\end{equation}
By the Stieltjes inversion formula~\cite{cit_10000_ST}, we obtain that $\mathbf{E}$ commutes with $B$.

\noindent
We denote by $\mathbf{M}(A,B)$ a set of all generalized resolvents $\mathbf{R}_z(A)$ of $A$ which satisfy
relation~(\ref{f3_43}).

Recall some known facts from~\cite{cit_4000_S} which we shall need here.
Let $K$ be a closed symmetric operator in a Hilbert space $\mathbf{H}$, with the domain $D(K)$,
$\overline{D(K)} = \mathbf{H}$.  Set $N_\lambda = N_\lambda(K) = \mathbf{H}
\ominus \Delta_K(\lambda)$, $\lambda\in \mathbb{C}\backslash \mathbb{R}$.

Consider an
arbitrary bounded linear operator $C$, which maps $N_i$ into $N_{-i}$.
For
\begin{equation}
\label{f3_45}
g = f + C\psi - \psi,\qquad f\in D(K),\ \psi\in N_i,
\end{equation}
we set
\begin{equation}
\label{f3_46}
K_C g = Kf + i C \psi + i \psi.
\end{equation}
Since an intersection of $D(K)$, $N_i$ and $N_{-i}$ consists only of the zero element,
this definition is correct.
Notice that $K_C$ is a part of the operator $K^*$.
The operator $K_C$ is said to be a {\it quasiself-adjoint extension of the operator $K$, defined by
the operator $K$}.

The following theorem can be found in~\cite[Theorem 7]{cit_4000_S}:
\begin{thm}
\label{t3_3}
Let $K$ be a closed symmetric operator in a Hilbert space $\mathbf{H}$ with the domain $D(K)$,
$\overline{D(K)} = \mathbf{H}$.
All generalized resolvents of the operator $K$ have the following form:
\begin{equation}
\label{f3_47}
\mathbf R_\lambda (K) = \left\{ \begin{array}{cc} (K_{F(\lambda)} - \lambda E_\mathbf{H})^{-1}, &
\mathop{\rm Im}\nolimits\lambda > 0\\
(K_{F^*(\overline{\lambda}) } - \lambda E_\mathbf{H})^{-1}, & \mathop{\rm Im}\nolimits\lambda < 0 \end{array}\right.,
\end{equation}
where $F(\lambda)$ is an analytic in $\mathbb{C}_+$ operator-valued function, which values are contractions
which map $N_i(A) = H_2$ into $N_{-i}(A) = H_4$ ($\| F(\lambda) \|\leq 1$),
and $K_{F(\lambda)}$ is the quasiself-adjoint extension of $K$ defined by $F(\lambda)$.

On the other hand, for any operator function $F(\lambda)$ having the above properties there corresponds by
relation~(\ref{f3_47}) a generalized resolvent of $K$.
\end{thm}
Notice that the correspondence between all generalized resolvents and functions $F(\lambda)$ in
Theorem~\ref{t3_3} is bijective~\cite{cit_4000_S}.

Return to the study of the Devinatz moment problem.
Let us describe the set $\mathbf{M}(A,B)$. Choose an arbitrary $\mathbf{R}_\lambda\in \mathbf{M}(A,B)$.
By~(\ref{f3_47}) we get
\begin{equation}
\label{f3_48}
\mathbf{R}_\lambda = (A_{F(\lambda)} - \lambda E_H)^{-1},\qquad \mathop{\rm Im}\nolimits\lambda > 0,
\end{equation}
where $F(\lambda)$ is an analytic in $\mathbb{C}_+$ operator-valued function, which values are contractions
which map $H_2$ into $H_4$, and
$A_{F(\lambda)}$ is the quasiself-adjoint extension of $A$ defined by $F(\lambda)$.
Then
$$ A_{F(\lambda)} = \mathbf{R}_\lambda^{-1} + \lambda E_H,\qquad \mathop{\rm Im}\nolimits\lambda > 0. $$
By virtue of relation~(\ref{f3_43}) we obtain
\begin{equation}
\label{f3_49}
BA_{F(\lambda)} h = A_{F(\lambda)} B h,\qquad h\in D(A_{F(\lambda)}),\ \lambda\in \mathbb{C}_+.
\end{equation}
Consider the following operators
\begin{equation}
\label{f3_50}
W_{\lambda} := (A_{F(\lambda)} + iE_H)(A_{F(\lambda)} - iE_H)^{-1} =
E_H + 2i(A_{F(\lambda)} - iE_H)^{-1},
\end{equation}
\begin{equation}
\label{f3_51}
V_A = (A +iE_H)(A - iE_H)^{-1} =
E_H + 2i(A - iE_H)^{-1},
\end{equation}
where $\lambda\in \mathbb{C}_+$.
Notice that (\cite{cit_4000_S})
\begin{equation}
\label{f3_52}
W_{\lambda} = V_A \oplus F(\lambda),\qquad \lambda\in \mathbb{C}_+.
\end{equation}
The operator $(A_{F(\lambda)} - iE_H)^{-1}$ is defined
on the whole $H$, see~\cite[p.79]{cit_4000_S}.
By relation~(\ref{f3_49}) we obtain
\begin{equation}
\label{f3_53}
B (A_{F(\lambda)} - iE_H)^{-1} h =
(A_{F(\lambda)} - iE_H)^{-1} B h,\qquad h\in H,\ \lambda\in \mathbb{C}_+.
\end{equation}
Then
\begin{equation}
\label{f3_54}
B W_\lambda = W_\lambda B,\qquad \lambda\in \mathbb{C}_+.
\end{equation}
Recall that by Proposition~\ref{p2_1} the operator $B$ reduces the subspaces $H_j$, $1\leq j\leq 4$,
and $BV_A = V_A B$. If we choose an arbitrary $h\in H_2$ and apply relations~(\ref{f3_54}),(\ref{f3_52}),
we get
\begin{equation}
\label{f3_55}
B F(\lambda) = F(\lambda) B,\qquad \lambda\in \mathbb{C}_+.
\end{equation}
Denote by $\mathbf{F}(A,B)$ a set of all analytic in $\mathbb{C}_+$ operator-valued functions
which values are contractions which map $H_2$ into $H_4$ and which satisfy relation~(\ref{f3_55}).
Thus, for an arbitrary $\mathbf{R}_\lambda\in \mathbf{M}(A,B)$ the corresponding function
$F(\lambda)\in \mathbf{F}(A,B)$. On the other hand, choose an arbitrary $F(\lambda)\in \mathbf{F}(A,B)$.
Then we derive~(\ref{f3_54}) with $W_\lambda$ defined by~(\ref{f3_50}). Then we get~(\ref{f3_53}),(\ref{f3_49})
and therefore
\begin{equation}
\label{f3_56}
B \mathbf{R}_\lambda  = \mathbf{R}_\lambda B,\qquad \lambda\in \mathbb{C}_+.
\end{equation}
Calculating the conjugate operators for the both sides of the last equality we conclude that this
relation holds for all $\lambda\in \mathbb{C}$.

\noindent
Consider the spectral measure $F_2$ of the operator $B_{H_2}$ in $H_2$. We have obtained
relation~(\ref{f3_39}) which we shall use one more time.
Notice that $F(\lambda)\in \mathbf{F}(A,B)$ if and only if
the operator-valued function
\begin{equation}
\label{f3_57}
G(\lambda) := W F(\lambda) U_{2,4}^{-1} W^{-1},\qquad \lambda\in \mathbb{C}_+,
\end{equation}
is analytic in $\mathbb{C}_+$ and has values which are
contractions in $\mathcal{H}$ which commute with $\mathbf{X}(\delta)$, $\forall\delta\in \mathfrak{[-\pi,\pi)}$.

This means that for an arbitrary $\lambda\in \mathbb{C}_+$ the operator
$G(\lambda)$ is decomposable and the values of the corresponding
operator-valued function $t(y)$ are $\mu$-a.e. contractions.
A set of all decomposable operators in $\mathcal{H}$ such that the values of the corresponding
operator-valued function $t(y)$ are $\mu$-a.e. contractions we denote by $\mathrm{T}(B;H_2)$.
A set of all analytic in $\mathbb{C}_+$ operator-valued functions $G(\lambda)$ with values
in $\mathrm{T}(B;H_2)$ we denote by $\mathbf{G}(A,B)$.

\begin{thm}
\label{t3_4}
Let a Devinatz moment problem~(\ref{f1_1}) be given. In conditions of Theorem~\ref{t3_1} all
solutions of the moment problem have the  form~(\ref{f3_4}) where the spectral
measures $\mathbf{E}$ of the operator $A$ are defined by the corresponding generalized
resolvents $\mathbf{R}_\lambda$ which are constructed by the following relation:
\begin{equation}
\label{f3_58}
\mathbf{R}_\lambda = (A_{F(\lambda)} - \lambda E_H)^{-1},\qquad \mathop{\rm Im}\nolimits\lambda > 0,
\end{equation}
where $F(\lambda) = W^{-1} G(\lambda) W U_{2,4}$, $G(\lambda)\in \mathbf{G}(A,B)$.

\noindent
Moreover, the correspondence between $\mathbf{G}(A,B)$ and a set of all solutions of
the Devinatz moment problem is bijective.
\end{thm}
{\bf Proof. } The proof follows from the previous considerations.
$\Box$

Consider an arbitrary non-negative Borel measure $\mu$ in the strip $\Pi$ which has all finite
moments~(\ref{f1_1}). What can be said about the density of power-trigonometric
polynomials~(\ref{f1_2}) in the corresponding space $L^2_\mu$?
The measure $\mu$ is a solution of the corresponding moment problem~(\ref{f1_1}).
Thus, $\mu$ admits a representation~(\ref{f3_4})
where $F$ is the spectral measure of $B$ and $\mathbf{E}$ is a spectral measure of $A$ which commutes with
$F$ (the  operators $A$ and $B$ in a Hilbert space $H$ are defined as above).

Suppose that (power-trigonometric) polynomials are dense in $L^2_\mu$.
Repeating arguments from the beginning of the Proof of Theorem~\ref{t3_1} we see that
in our case $H_0 = \{ 0 \}$ and $\widetilde A$, $\widetilde B$ are operators in $H$.
Moreover, we have $\mu = ((\widetilde E\times \widetilde F) x_{0,0}, x_{0,0})_{H}$,
where $\widetilde E$ is the spectral measure of $\widetilde A$, $\widetilde F = F$.
Consequently, $\mu$ is  a canonical solution of the Devinatz moment problem.

\noindent
The converse assertion is more complicated and will be studied elsewhere.

\vspace{1.5cm}

Sergey M. Zagorodnyuk

School of Mathematics and Mekhanics

Karazin Kharkiv National University

Kharkiv, 61077

Ukraine

\begin{center}
\bf
Devinatz's moment problem: a description of all solutions.
\end{center}

\begin{center}
\bf
S.M. Zagorodnyuk
\end{center}

In this paper we study Devinatz's moment problem:
to find a non-negative Borel measure $\mu$ in a strip
$\Pi = \{ (x,\varphi):\ x\in \mathbb{R},\ -\pi\leq \varphi < \pi \},$
such that $\int_\Pi x^m e^{in\varphi} d\mu = s_{m,n}$, $m\in \mathbb{Z}_+$, $n\in \mathbb{Z}$,
where $\{ s_{m,n} \}_{m\in \mathbb{Z}_+, n\in \mathbb{Z}}$ is a given sequence of complex numbers.
We present a  new proof of the Devinatz solvability criterion for this moment problem.
We obtained a parameterization of all solutions of Devinatz's moment problem.
We used an abstract operator approach and results of Godi\v{c}, Lucenko and
Shtraus.

\vspace{1cm}
Key words: moment problem, measure, generalized resolvent.
\vspace{1cm}

MSC 2000: 44A60, 30E05.

}
\end{document}